\documentclass[12pt]{article}
\usepackage{graphicx, tikz, color,url}
\usepackage{amsmath,amssymb,latexsym,amsthm}
\usepackage{amsfonts}
\usetikzlibrary{arrows}
\usetikzlibrary{positioning, quotes}
\usepackage{soul}

\newtheorem{theorem}{Theorem}[section]

\newtheorem{lemma}[theorem]{Lemma}
\newtheorem{proposition}[theorem]{Proposition}
\newtheorem{observation}[theorem]{Observation}

\newcommand{\gsmb}{\gamma_{\rm SMB}}
\newcommand{\gmb}{\gamma_{\rm MB}}
\newcommand{\smallqed}{{\tiny ($\Box$)}}

\begin{document}
\title{On Maker-Breaker domination game critical graphs}
\author{
Boštjan Brešar$^{a,b}$
\and
Tanja Dravec$^{a,b}$
\and
Kirsti Kuenzel$^{c}$\\
\and
Douglas F.\ Rall$^{d}$
}
\maketitle

\begin{center}
$^a$ Faculty of Natural Sciences and Mathematics, University of Maribor, Slovenia\\
$^b$ Institute of Mathematics, Physics and Mechanics, Ljubljana, Slovenia\\
$^c$ Department of Mathematics, Trinity College, Hartford, CT, USA\\
$^d$ Emeritus Professor of Mathematics, Furman University, Greenville, SC, USA\\
\end{center}
\medskip

\begin{abstract}
The Maker-Breaker domination game is played on a graph $G$ by  Dominator and Staller who alternate turns selecting an unplayed vertex of $G$. The goal of Dominator is that the vertices he selected during the game form a dominating set while Staller's goal is to prevent this from happening. The graph invariant $\gmb'(G)$ is the number of Dominator's moves in the game played on $G$ in which he can achieve his goal  when Staller makes the first move and both players play optimally. In this paper, we continue the investigation of $2$-$\gmb'$-critical graphs, initiated in [Divarakan et al., Maker--Breaker domination game critical graphs, Discrete Appl.\  Math. 368 (2025) 126--134], which are defined as the graphs $G$ with $\gmb'(G)=2$ and $\gmb'(G-e)>2$ for every edge $e$ in $G$. The authors characterized  bipartite $2$-$\gmb'$-critical graphs, and found an example of a non-bipartite $2$-$\gmb'$-critical graph. In this paper, we characterize the $2$-$\gmb'$-critical graphs that have a cut-vertex, which are represented by two infinite families. In addition, we
prove that $C_5$ is the only non-bipartite, triangle-free $2$-$\gmb'$-critical graph.
\end{abstract}


\noindent
{\bf AMS Subj.\ Class.\ (2020)}: 05C57, 05C69

 \section{Introduction}
\label{sec: intro}

Maker-Breaker games are an important and ever-increasing area of study in positional games.  Many of these games are played on a graph or a hypergraph;  see a monograph on positional games~\cite{hefetz-2014} and a recent study of its computational aspects~\cite{rahman-2023}.  
Other games played on graphs or hypergraphs have been introduced and studied.  The recent book~\cite{book-2021} concerns two ``domination type'' games played on graphs.
The \emph{Maker-Breaker domination game} was introduced in 2020 by Duch\^{e}ne, Gledel, Parreau, and Renault~\cite{duchene-2020} as follows.  Let $G$ be a finite, simple graph.  Two players, Dominator and Staller take turns selecting a vertex of $G$ that has not been previously chosen in the game.  
Dominator's goal is that at some point in the course of play, the set of vertices he has selected is a dominating set of $G$. Staller's goal is to prevent Dominator from constructing such a dominating set.  The game ends when one of the two players accomplishes the goal. (Recall that $D\subseteq V(G)$ is a {\em dominating set} of $G$ if every vertex in $V(G)-D$ is adjacent to a vertex in $D$, and $\gamma(G)$ is the cardinality of a smallest such set.) There are two versions of the game. If Dominator (Staller) is the first player to select a vertex, then we call this the \emph{D-game (S-game)}.  

The Maker-Breaker domination game attracted a number of researchers studying it from various perspectives~\cite{bagan2025partition, bujtas-2023,divakaran-2024,dokyeesun-2024+,forcan-2023}. In addition, a version of the game in the context of total domination was introduced in~\cite{gledel-2020}. 
The authors of~\cite{gledel-2019} introduced four graph invariants that arise naturally from these games (see also~\cite{bujtas-2024, forcan-2023} for their further studies).  In this context, the goal of the player who wins the game is to minimize the number of vertices he or she selects during the entire game while the goal of the other player is to maximize that number.  If Dominator wins the D-game (S-game) played on $G$, then the number of moves he makes when both players are playing optimally is the {\em Maker-Breaker domination number}, which is denoted by $\gmb(G)$ (respectively, $\gmb'(G)$ when Staller starts the game). We set $\gmb(G)=\infty$ (respectively, $\gmb'(G)=\infty$ in the S-game) if Staller wins the game on $G$. In a similar way, counting the number of moves Staller needs to win the D-game (S-game) the invariants denoted by $\gsmb(G)$ (respectively, $\gsmb'(G)$) arise. 

Recently,  Divakaran, Dravec, James, Klav{\v{z}}ar and Nair~\cite{ddj-2025} introduced criticality concepts with respect to Maker-Breaker domination numbers. 
Notably, a graph $G$ is {\em $\gmb$-critical} (respectively, {\em $\gmb'$-critical}) if  $\gmb(G-e)>\gmb(G)$ (respectively, $\gmb'(G-e)>\gmb'(G)$) holds for every $e\in E(G)$. If $\gmb(G)=k$ and $G$ is $\gmb$-critical, then $G$ is {\em $k$-$\gmb$-critical}. The {\em $k$-$\gmb'$-critical} graphs are defined analogously. 
We will focus solely on these two classes of graphs. For the definitions of $\gsmb$- and $\gsmb'$-critical graphs see~\cite{ddj-2025}. 
In~\cite{ddj-2025} the existence of $k$-$\gmb$-critical graphs and $k$-$\gmb'$-critical graphs is established for all $k\ge 2$. In addition, critical graphs are characterized for most of the cases in which Dominator wins the game in one or two moves. The main exception are  $2$-$\gmb'$-critical graphs for which only a partial characterization is provided in~\cite{ddj-2025}, i.e.\  for the case when $G$ is bipartite. The authors also found a non-bipartite $2$-$\gmb'$-critical graph and posed the problem~\cite[Problem 4.8]{ddj-2025} of characterizing this class of graphs. In this paper, we continue the study on $2$-$\gmb'$-critical graphs. We find an infinite family of non-bipartite $2$-$\gmb'$-critical graphs, characterize the $2$-$\gmb'$-critical graphs that have a cut-vertex, and prove that the only non-bipartite triangle-free $2$-$\gmb'$-critical graph is $C_5$.

\section{Definitions and preliminaries}

Let $G$ be a  graph. The order of  $G$ is denoted by $n(G)$. For a vertex $v\in V(G)$,  its neighborhood is denoted by $N(v)$ and its closed neighborhood by $N[v]$.  The degree of $v$  is $\deg(v)=|N(v)|$.  The minimum and the maximum degree of $G$ are respectively denoted by $\delta(G)$ and $\Delta(G)$. An {\em isolated vertex} is a vertex of degree $0$, a {\em leaf} is a vertex of degree $1$.  For a graph $G$ and $S \subseteq  V(G)$ let $G[S]$ denote the subgraph of $G$ induced by $S$. In addition, $G-S$ denotes $G[V(G)-S]$. For $S=\{x\}$, $G-\{x\}$ is shortly denoted by $G-x$. For $F \subseteq E(G)$, $G-F$ denotes the spanning subgraph of $G$ obtained from $G$ by deleting edges in $F$. For $F=\{e\}$ we write $G-e$ instead of $G-\{e\}$.
A vertex $x$ of a connected graph $G$ is a {\em cut-vertex} if $G-x$ is not connected. An edge $e$ of a connected graph $G$ is a {\em bridge} if $G-e$ is not connected. If $G$ and $H$ are graphs, then the {\em join} of $G$ and $H$, denoted by $G \vee H$, is constructed from their disjoint union by adding all edges of the form $uv$ where $u\in V(G)$ and $v \in V(H)$.

A {\em dominating set} of $G$ is a set $D \subseteq V(G)$ such that each vertex from $V(G)-D$ has a neighbor in $D$. The {\em domination number} $\gamma(G)$ of $G$ is the minimum cardinality of a dominating set of $G$. A vertex of degree $n(G) -1$ is a {\em universal vertex} of $G$, and the set of universal vertices in $G$ is denoted by $\mathfrak{U}(G)$. The cardinality of $\mathfrak{U}(G)$ is denoted by $\mathfrak{u}(G)$.  We will denote the sequence of vertices chosen when the S-game is played on a graph by $s_1,d_1,s_2,d_2, \ldots$, where $s_i$ ($d_i$) is the $i^{{\rm th}}$ vertex played by Staller (respectively, by Dominator).

We recall known results needed in the remainder. For $X\subseteq V(G)$, let $G|X$ denote the graph $G$ in which vertices from $X$ are considered as being already dominated. Then we have: 

\begin{theorem} [Continuation Principle~\cite{gledel-2019}]
\label{thm: Continuation Principle}
Let $G$ be a graph with $A, B \subseteq V(G)$. If $B\subseteq A$ then $\gmb(G|A) \leq \gmb(G|B)$ and  $\gmb'(G|A) \leq \gmb'(G|B)$.
\end{theorem}

\begin{lemma}{\rm \cite{dokyeesun-2024+}}\label{lem:doky}
If $G$ is a graph, then the following properties hold: 
 \begin{enumerate}
 \item[(i)] $\gmb(G)\leq \gmb(G-e)$ for every $e\in E(G)$,
\item[(ii)] $\gmb'(G)\leq \gmb'(G-e)$ for every $e \in E(G)$.
\end{enumerate}
\end{lemma}

\subsection*{Families of connected $2$-$\gmb'$-critical  graphs}

 In~\cite{ddj-2025} the following necessary conditions were proved for connected $2$-$\gmb'$-critical graphs.

\begin{proposition} \label{Prop:necessary-conditions}\cite{ddj-2025}
If $G$ is a connected $2$-$\gmb'$-critical graph, then the following properties hold. 
\begin{enumerate}
    \item[(i)] $n(G)\ge 5$.
    \item[(ii)] $\delta(G) \geq 2$.
    \item[(iii)] $\Delta(G) \leq n(G)-2$.
\end{enumerate}
\end{proposition}

If $G$ is a connected $2$-$\gmb'$-critical graph, then Proposition~\ref{Prop:necessary-conditions}(iii) implies that $\gamma(G) \geq 2$. Hence Dominator cannot win the S-game in one move, yet regardless of Staller's moves he wins in two moves. The following observation follows.

\begin{observation}\label{o:1}
If $G$ is a connected $2$-$\gmb'$-critical graph, then any unplayed vertex in $G$ is an optimal move for Staller in the S-game.
\end{observation}

Connected bipartite $2$-$\gmb'$-critical graphs were also characterized in~\cite{ddj-2025}. Let us introduce the family $\cal B$ as follows. A graph $G \in \cal B$ if $G=K_{2,m}$, $m\geq 3$, or $G$ is a bipartite graph with a bipartition $V(G)=V_1 \cup V_2$, where $|V_1|=m \geq 3$, $|V_2|=n \geq 3$, and $V_1$ contains exactly two vertices of degree $n$, $V_2$ contains exactly two vertices of degree $m$, while all the other vertices are of degree 2.

\begin{theorem}\label{thm:bip}~\cite{ddj-2025}
Let $G$ be a connected bipartite graph. Then $G$ is $2$-$\gmb'$-critical if and only if $G\in \cal B$.
\end{theorem}

\begin{figure}
\begin{center}

\begin{tikzpicture}[scale=0.85,main_node/.style={circle,draw,minimum size=1em,inner sep=5pt]}]

\draw (-3.5, 0) node {$v_1$};
\draw (-0.5,0) node {$v_2$};

\node[main_node] (0) at (-5, 1.6) {};
\node[main_node] (1) at (-5, -1.6) {};
\node[main_node] (2) at (-3.5, 0) {};
\node[main_node] (3) at (-2,0) {};
\node[main_node] (4) at (-0.5,0) {};
\node[main_node] (5) at (1, 1.6) {};
\node[main_node] (6) at (1,-1.6) {};


\draw (0)--(1)--(2)--(0);
\draw (2)--(3)--(4)--(5)--(6)--(4);

\draw (5, 0) node {$v_1$};
\draw (8,0) node {$v_2$};

\node[main_node] (a) at (5, 0) {};
\node[main_node] (b1) at (6.5,1.6) {};
\node[main_node] (b2) at (6.5, 0) {};
\node[main_node] (b3) at (6.5,-1.6) {};
\node[main_node] (c) at (8,0) {};
\node[main_node] (d1) at (9.5, 1.6) {};
\node[main_node] (d2) at (9.5,-1.6) {};
\node[main_node] (e) at (11,0) {};

\draw (a)--(b1)--(c)--(b2)--(a)--(b3)--(c);
\draw (c)--(d1)--(d2)--(c)--(d2)--(e)--(d1);
\draw (e)--(c);

\end{tikzpicture}
\end{center}
\caption{Graph $F_{2,1,2}$ and graph $F_{0,3,3}$.}
    \label{fig:twoexamples}
\end{figure}
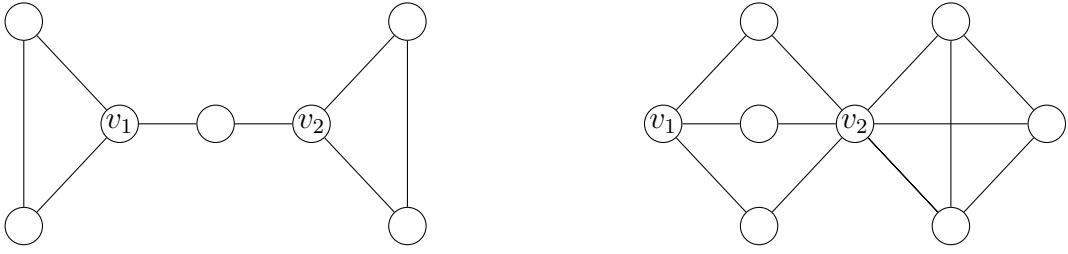

Finally, the authors in~\cite{ddj-2025} surprisingly noted that there also exist connected, non-bipartite $2$-$\gmb'$-critical graphs and presented one example of such a graph, which is the graph $G$ obtained from two disjoint copies of $K_3$ and a vertex $x$ by connecting $x$ with exactly one vertex in each copy of $K_3$ (see the left graph in Fig.~\ref{fig:twoexamples}). Next, we present an infinite family of non-bipartite $2$-$\gmb'$-critical graphs that contains $G$.

Consider the family of graphs $\mathcal{F}$ defined as follows. Let $H_m = K_2 \vee \overline{K_{m-2}}$ if $m\ge 2$ (where in the case $m=2$ we mean $H_2=K_2$), and if $m=0$, $H_0$ is the null graph. (That is, $H_0$ is the graph with no vertices and no edges.)


\begin{figure}[h!]
\begin{center}
    
\begin{tikzpicture}[scale=1.2, main_node/.style={circle,draw,minimum size=1em,inner sep=5pt]}]

\draw (0,5) node {}; 
\draw (-3.742857251848493, 3.2) node {$v_1$};
\draw (-1.0285715375627793, 3.2) node {$v_2$};
\node[main_node] (0) at (-4.9, 3.7571430751255583) {};
\node[main_node] (1) at (-4.9, 2.242857360839844) {};
\node[main_node] (2) at (-3.742857251848493, 3.2) {};
\node[main_node] (3) at (-2.4, 4.4) {};
\node[main_node] (4) at (-2.4, 3.2) {};
\node[main_node] (5) at (-2.4, 2) {};
\node[main_node] (6) at (-1.0285715375627793, 3.2) {};
\node[main_node] (7) at (-0.11428582327706494, 4.071428571428571) {};
\node[main_node] (8) at (-0.05714296613420755, 2.1571430751255587) {};
\node[main_node] (9) at (1.1, 4.414285714285715) {};
\node[main_node] (10) at (1.1, 3.2142859322684147) {};
\node[main_node] (11) at (1.1, 1.8285716465541295) {};


 \path[draw, thick]
(0) edge node {} (1) 
(1) edge node {} (2) 
(2) edge node {} (0) 
(2) edge node {} (3) 
(2) edge node {} (4) 
(2) edge node {} (5) 
(3) edge node {} (6) 
(6) edge node {} (4) 
(6) edge node {} (5) 
(6) edge node {} (7) 
(6) edge node {} (8) 
(7) edge node {} (8) 
(8) edge node {} (11) 
(11) edge node {} (7) 
(9) edge node {} (7) 
(7) edge node {} (10) 
(8) edge node {} (10) 
(8) edge node {} (9) 
(6) edge node {} (9) 
(6) edge node {} (10) 
(6) edge node {} (11) 
;

\end{tikzpicture}
\end{center}
\caption{Graph $F_{2,3,5}$.}
    \label{fig:family}
\end{figure}
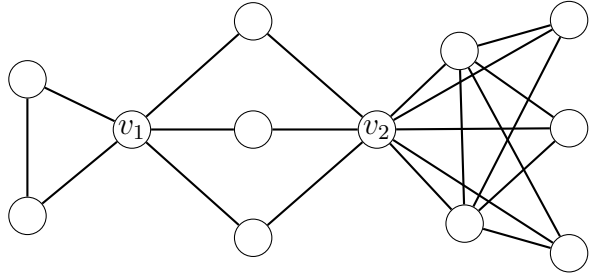

Let $t \in \mathbb{N}$ and let 
$\{m_1, m_2\} \subset (\mathbb{N}- \{1\})\cup \{0\}$, such that 
$m_1 \le m_2$ and $m_2\geq 2$, and
let $\{v_1, v_2\}$ and $B$ be independent sets such that $|B|=t$. In addition, if $m_1=0$, then $t\ge 2$.   
The graph $F_{m_1,t,m_2}$ is obtained from the disjoint union of  $H_{m_1}$, $H_{m_2}$, $B$, and 
$\{v_1,v_2\}$ by adding an edge between $v_1$ and each vertex in $H_{m_1}$, adding an edge between $v_2$ and each vertex in $H_{m_2}$, and adding all edges between $\{v_1, v_2\}$ and $B$.  See Figs.~\ref{fig:twoexamples} and~\ref{fig:family}. Note that $F_{2,1,2}$ was found by the authors in~\cite{ddj-2025}. A graph $G$ is in the family $\mathcal{F}$ if and only if $G=F_{m_1,t,m_2}$ for some choice of $m_1,t,$ and $m_2$ as specified above. 

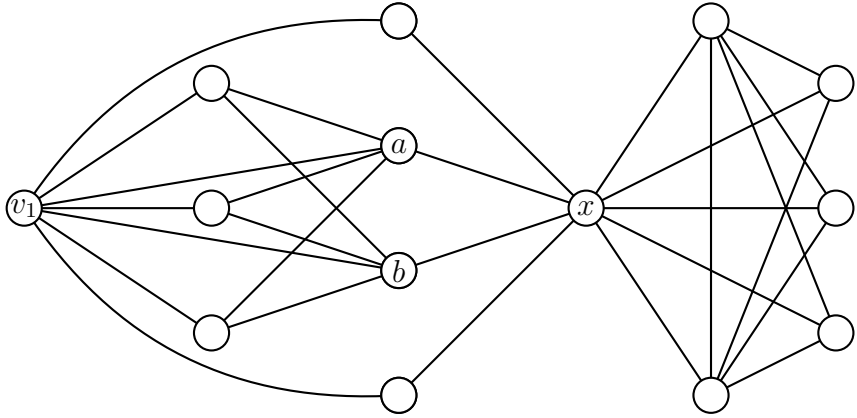
\begin{figure}[ht!]
\begin{center}
\begin{tikzpicture}[scale=0.83,style=thick]
\def\vr{8pt} 
\draw (-7,0) to [bend left=30] (-1,3);
\draw (-7,0) to [bend right=30] (-1,-3);
\path (-7,0) coordinate (v1);
\path (-4,2) coordinate (u1); \path (-4,0) coordinate (u2);  \path (-4,-2) coordinate (u3);
\path (-1,3) coordinate (w1); \path (-1,-3) coordinate (w2); \path (-1,1) coordinate (w3); \path (-1,-1) coordinate (w4);
\path (2,0) coordinate (x); \path (4,3) coordinate (h1); \path (4,-3) coordinate (h2);
\path (6,2) coordinate (y1); \path (6,0) coordinate (y2); \path (6,-2) coordinate (y3);

\foreach \i in {1,...,3}
{  \draw (u\i) -- (v1); \draw (u\i) -- (w3); \draw (u\i) -- (w4); }
\foreach \i in {3,4}
{  \draw (w\i) -- (v1); }
\foreach \i in {1,...,4}
{  \draw (w\i) -- (x); }
\foreach \i in {1,...,2}
{  \draw (h\i) -- (x); }
\foreach \i in {1,...,3}
{  \draw (y\i) -- (x); }
\foreach \i in {1,...,3}
{  \draw (y\i) -- (h1); }
\foreach \i in {1,...,3}
{  \draw (y\i) -- (h2); }
\draw (h1) -- (h2);

\foreach \i in {1,...,3}
{  \draw (u\i)  [fill=white] circle (\vr); \draw (y\i)  [fill=white] circle (\vr); }
\foreach \i in {1,...,2}
{  \draw (w\i)  [fill=white] circle (\vr); }
\foreach \i in {1,...,4}
{  \draw (w\i)  [fill=white] circle (\vr); }
\draw (v1)  [fill=white] circle (8pt); \draw (x)  [fill=white] circle (8pt); \draw (h1)  [fill=white] circle (\vr); \draw (h2)  [fill=white] circle (\vr);
\draw (w4)  [fill=white] circle (8pt); \draw (w3)  [fill=white] circle (8pt);
\draw (-7,0) node {$v_1$}; \draw (2,0) node {$x$}; \draw (-1,1) node {$a$}; \draw (-1,-1) node {$b$};

\end{tikzpicture}
\end{center}
\caption{Graph $F'_{3,4,5}$.}
    \label{fig:familyprime}
\end{figure}

Next the family of graphs $\mathcal{F}'$ is defined as follows. 
Let $\{s,q,m\} \subset \mathbb{N}-\{1\}$. 
Let graph $P$ be obtained from the graph $K_{2,s}$, where $\{a,b\}$ is a part of the bipartition of $K_{2,s}$, by adding $q-2$ vertices (if $q>2$) $b_1,\ldots, b_{q-2}$, which induce $\overline{K_{q-2}}$. Now, let $Q$ be the graph $K_1\vee P$, where $v_1$ is the vertex representing $K_1$ in the definition of $Q$.  The graph $F'_{s,q,m}$ is obtained from the disjoint union of $Q$ and the graph $H_m$ by adding a  vertex $x$ and connecting it by an edge to each vertex in $V(H_m)\cup\{a,b,b_1,\ldots,b_{q-2}\}$. 
See Fig.~\ref{fig:familyprime} depicting $F'_{3,4,5}$. A graph $G$ is in the family $\mathcal{F}'$ if and only if $G=F'_{s,q,m}$ for some choice of $s,q,$ and $m$ as specified above.
\textbf{}

Now that the families $\cal F$ and ${\cal F}'$ are defined, we can state the main results of the paper. We will give their proofs in Section~\ref{sec:triangle-free} and Section~\ref{sec:proofC5}.

\begin{theorem}\label{thm:NoK3}
    If $G$ is a connected, triangle-free graph, then $G$ is $2$-$\gmb'$-critical graph if and only if $G \in \cal B$ or $G = C_5$. 
\end{theorem}

\begin{theorem}
\label{thm:cutvertex}
    If $G$ is a connected graph with a cut-vertex, then $G$ is $2$-$\gmb'$-critical if and only if $G \in \mathcal{F}\cup\mathcal{F}'$.
\end{theorem}

We complete this section with an auxiliary result that will be used several times in our proofs. Note that if we omit the condition that $G$ is connected in Proposition~\ref{Prop:necessary-conditions}, then item (ii) in that proposition does not hold true. Notably, there exist $2$-$\gmb'$-critical graphs $G$ with $\delta(G)=1$, which we characterize in the next lemma. 
 In the proof we use the result from~\cite{ddj-2025} that $1$-$\gmb'$-critical graphs are exactly the graphs $H_m$ for $m\geq 2$.

\begin{lemma}\label{l: NoLeaves}
    Let $G$ be a graph with $\delta(G)=1$. Then $G$ is $2$-$\gmb'$-critical if and only if $G$ is the disjoint union of $H_2$ and $H_{m}$ for some $m\geq 2$. 
\end{lemma}
\begin{proof}
    Let $G$ be a $2$-$\gmb'$-critical graph with $\delta(G)=1$. Then it follows from Proposition~\ref{Prop:necessary-conditions} that $G$ is not connected. Since $\gmb'(G)=2$, $G$ has two components each having at least 2 universal vertices. Since each component must be $1$-$\gmb'$-critical  and $\delta(G)=1$, one component must be $H_2$ and the other component is isomorphic to $H_m$ for some $m \geq 2$.  The proof of the converse is straightforward and is left to the reader.
\end{proof}

\section{Proof of Theorem~\ref{thm:NoK3}}
\label{sec:triangle-free}

In this section we prove Theorem~\ref{thm:NoK3}.

Let $G$ be a connected, triangle-free 
$2$-$\gmb'$-critical graph. Let $s_1, d_1, s_2, d_2, \ldots$ be a sequence of optimal moves  and suppose that $d_1d_2 \in E(G)$. This implies that $V(G)=N_G(d_1) \cup N_G(d_2)$ and $N_G(d_1) \cap N_G(d_2) = \emptyset$, as $G$ is triangle-free. Therefore, we can partition $V(G)$ in two sets $A = N_G(d_1)$ and $B = N_G(d_2)$, both being independent sets. Hence $G$ is bipartite and thus by Theorem~\ref{thm:bip} $G \in \cal B$.

Therefore, we may assume that there is no sequence of optimal moves such that $d_1d_2 \in E(G)$. After Dominator plays $d_1$, it must be that both $s_2$ and $d_2$ dominate $V(G) - N_G[d_1]$ 
for otherwise $d_2$ would be the optimal play for Staller on her second move. Since $d_2 \in V(G) - N_G[d_1]$, this means that $s_2d_2 \in E(G)$. This implies that $|V(G) - N_G[d_1]| \le 2$ for otherwise there is some $w \in V(G) - (N_G[d_1]\cup \{s_2, d_2\})$ and then $ws_2d_2$ is a triangle in $G$. If $|V(G)-N_G[d_1]|=1$, then $V(G)-N_G[d_1]=\{d_2\}$.  In this case $d_2$ dominates $N_G(d_1)$, for otherwise some vertex in $N_G(d_1)$ has degree $1$, which contradicts Proposition~\ref{Prop:necessary-conditions}. Thus $G$ is bipartite and so $G \in \cal B$. Now, we may assume that $|V(G) - N_G[d_1]| = 2$. We may write $V(G) - N_G[d_1] = \{d_2, a\}$, and note that $d_2a\in E(G)$. Since $G$ is triangle-free, $N_G(d_1)$ is an independent set and we may partition $N_G(d_1) = A \cup B \cup C$ where $A = N_G(a) \cap N_G(d_1)$, $B = N_G(d_1) \cap N_G(d_2)$ and $C = N_G(d_1) - (N_G(a) \cup N_G(d_2))$. 
We claim that $C = \emptyset$. Indeed, suppose $w \in C$. Since $w$ is not adjacent to any vertex in $\{a, d_2\}$, it has a neighbor $t \in N_G(d_1) - \{w\}$ as $\delta(G) \ge 2$ by Proposition~\ref{Prop:necessary-conditions} and Lemma~\ref{l: NoLeaves}. However, this implies $d_1tw$ is a triangle. Hence, $C = \emptyset$. If $|A| = 1 = |B|$, then $G = C_5$. 

Without loss of generality, assume that $|A| \ge 2$. Consider the game when Staller's first move is $d_1$. In order to dominate $d_1$, Dominator must choose a vertex from $A\cup B$ on either his first or second move. If he chooses a vertex from $A \cup B$ on his first move, say $w$, then Staller replies by choosing $a$, and Dominator needs at least two more moves to dominate all vertices of $G$, a contradiction.
Thus, Dominator chooses a vertex from $\{a, d_2\}$ in his first move, and a vertex from $A \cup B$ in his second move. If he chooses $d_2$ in his first move, then $G - N_G[d_2]$ is a star $K_{1, |A|}$ with center $d_1$, which Dominator cannot play to win the game. Therefore, Dominator must choose $a$ in his first move, and the only possibility for his win in two moves is that $|B| = 1$; say $B = \{b\}$. Now, Staller can choose $b$ in her second move in which case Dominator cannot win in his second move. Indeed, Dominator must dominate both $d_1$ and $b$, and yet he cannot choose either $d_1$ or $b$. Thus, $\gmb'(G)>2$, the final contradiction. 

The converse follows from Theorem~\ref{thm:bip} and the fact that $C_5$ is $2$-$\gmb'$-critical. 

\section{Proof of Theorem~\ref{thm:cutvertex}}\label{sec:proofC5}

This section is devoted to the proof of Theorem~\ref{thm:cutvertex}. We will present a series of auxiliary results, which will be combined into a straightforward proof of the theorem at the end of this section.

First, let us prove the sufficient condition in Theorem~\ref{thm:cutvertex}.

\begin{proposition}\label{prp:familyF} If $G \in \mathcal{F}$, then $G$ is a $2$-$\gmb'$-critical graph.
\end{proposition}
\begin{proof}
    Let $G \in \cal{F}$. Then $G=F_{m_1,t,m_2}$ for some nonnegative integers $m_1,t,m_2$ satisfying the conditions given in the definition of the family $\mathcal{F}$.
    Since $\{v_1,v_2\}$ dominates $G$ and $\mathfrak{u}(G)=0$, we have $\gamma(G)=2$.
    
    Now we present Dominator's strategy to win in two moves.
    Dominator in his first move selects vertex $v_1$ if it is possible, otherwise he selects $v_2$. If the first case is possible, then in his second move he selects one among at least three universal vertices of the subgraph $G[\{v_2\} \cup V(H_{m_2})]$, which ends the game. Now, suppose Dominator is forced to select $v_2$ in his first move.  If $m_1=0$, then by definition $t\ge 2$, and so in his second move Dominator can select a vertex $b\in B$, and $\{v_2,b\}$ dominates $G$. If $m_1 \geq 2$, then in his second move he selects one of the universal vertices of $H_{m_1}$ and finishes the game. Thus, in all cases $\gmb'(G) =2$.  

    It remains to show that $G$ is critical. If $e$ is incident with a vertex from $B$, then $\delta(G-e)=1$. Suppose $\gmb'(G-e)=2$. Then $G-e$ contains a spanning $2$-$\gmb'$-critical subgraph $G'$ with minimum degree $1$. (Clearly, $\delta(G')=0$ is impossible.).  Thus, by Lemma~\ref{l: NoLeaves}, 
 $G'$ is the disjoint union of two graphs $H_2$ and $H_b$ for some $b \geq 2$. Such a graph cannot be a spanning subgraph of $G$, which is a contradiction. Hence $\gmb'(G-e) \geq 3$. 

 Now let $e$ be an edge connecting $v_1$ to $H_{m_1}$ or $v_2$ to $H_{m_2}$ and denote $G-e$ by $G'$.  
 Without loss of generality we may assume that $e$ is incident with $v_2$. Let $\{u_1,u_2\}\subseteq \mathfrak{U}(H_{m_2})$.  Suppose first that $e$ is incident with a vertex in $\mathfrak{U}(H_{m_2})$, say $e=v_2u_1$.  
 Staller's strategy is to select $u_2$ as her first move. If $m_2=3$, then $H_{m_2} = K_3$ and let $ \mathfrak{U}(H_{m_2})=\{u_1,u_2,u_3\}= V(H_{m_2})$.  In this case Dominator would need to play both vertices of the set $\{u_3,v_1\}$ to finish the game in two moves.  Staller's second move, then, is to select the vertex remaining in  $\{u_3,v_1\}$ after Dominator's first move, which shows in this case that $\gmb'(G') > 2$.  If $m_2 \neq 3$, then $\gamma((G-e)-u_2) > 2$ and it follows that $\gmb'(G') > 2$ as desired. 
  Assume now that $e$ is an edge between $v_2$ and $H_{m_2}$ not incident with a vertex in $\mathfrak{U}(H_{m_2})$. Denote $e=v_2v$. Staller's strategy is to select a vertex from $\{u_1,u_2\}$, say $u_1$, as her first move. If Dominator does not select $u_2$ in his first move, then Staller can choose $u_2$ in her second move so that Dominator cannot dominate $\{v_2\} \cup V(H_{m_2})$ with one vertex. If Dominator selects $u_2$ in his first move, then Staller can select $v_1$ in her second move. In this case,  Dominator cannot dominate $V(H_{m_1})\cup B$ with one vertex even when $m_1=0$, it must be that $|B| \ge 2$. Hence, $\gmb'(G') \geq 3$.



Finally assume that $e$ is an edge of $H_{m_1}$ or $H_{m_2}$, say $H_{m_2}$. 
Hence $e$ is incident to a universal vertex, say $u_2$, of $H_{m_2}$.  In this case, Staller selects $v_2$ in her first move. First, assume that $m_2\ne 3$. Then, in her second move Staller selects a vertex in $\{u_1,v_1\}$ depending on the first choice of Dominator. The reader can easily verify that the presented strategy of Staller leads to at least three moves of Dominator. Now, let $m_2=3$.  Then, in her second move Staller selects a vertex in $\{u,v_1\}$ depending on the first choice of Dominator, where $u$ is the only universal vertex of $H_{m_2}-e$. The argument that three moves of Dominator are required is the same as the above. 
 \end{proof}

\begin{proposition}
\label{prp:familyFprime} If $G \in \mathcal{F}'$, then $G$ is a $2$-$\gmb'$-critical graph.
\end{proposition}
\begin{proof}
Let $G \in \mathcal{F}'$. By definition, $G=F'_{s,q,m}$ for some choice of the positive integers $s,q,$ and $m$ as specified in the description of $\mathcal{F}'$.  Let $\{h_1,h_2\} \subseteq \mathfrak{U}(H_m)$. If $m=3$, then the third vertex of $H_m$, say $h_3$, is also a universal vertex of $H_m$.   If $m=2$ or $m>3$, then it is straightforward to show that 
\[\{x,v_1\}, \{x,a\}, \{x,b\}, \{v_1,h_1\}, \{v_1,h_2\}\] is a complete list of the minimum dominating sets of $G$.  For $m=3$, adding the set $\{v_1,h_3\}$ is then such a complete list of the $\gamma$-sets of $G$.

We first show that $\gmb'(G)=2$ by specifying a strategy for Dominator that will end the game after his second move.  If Staller plays $x$ on her first move, then Dominator plays $d_1=v_1$. Regardless of Staller's second move, Dominator can play $d_2 \in \{h_1,h_2\}$ to finish the game.  Similarly, if $s_1=v_1$, then Dominator will play $d_1=x$ and either $a$ or $b$ as his second move.  If $s_1 \in \{a,b\}$, then Dominator will play $d_1=v_1$ and then $x,h_1$ or $h_2$ on his second move.  If $s_1 \in \mathfrak{U}(H_m)$, then Dominator will play $d_1=x$ followed by $v_1,a$ or $b$ on his second move.  Finally, if $s_1 \notin \{x,v_1,a,b\} \cup \mathfrak{U}(H_m)$, then Dominator will play $v_1$ on his first move.  In this case, by selecting $d_2\in \{x,h_1,h_2\}$ Dominator will finish the game in two moves.  Therefore, $\gmb'(G)=2$.

Now we show that $G$ is critical.  That is, we will prove that $\gmb'(G-e) > 2$ for every $e \in E(G)$.  Assume first that $m \neq 3$. For each edge $e$, we will show that the dominating sets of cardinality $2$ for $G-e$ is a proper subset of the original list of $\gamma$-sets of $G$.  This will allow a choice for $s_1$ so that Dominator cannot finish the game in two moves.  

Let us now consider the different possibilities for the edge $e$ that will be removed from $G$.  Assume first that $e=xw$ for some $w \in V(H_m) \cup \{b_1,\ldots,b_{q-2}\}$. The vertex  $v_1$ belongs to every $\gamma$-set of $G-e$.  In this case Staller will select $v_1$ as her first move, which means that Dominator cannot finish the game with his first two moves.  
Next, let $e$ be incident with $v_1$.  Staller can select $x$ on her first move.  Since $x$ belongs to every $\gamma$-set of $G-e$, Dominator will require at least three moves to dominate $G-e$.  
Now let $e=h_1h_2$. Since $x$ belongs to all $\gamma$-sets of $G-e$,  Staller plays $s_1=x$.  Thus, Dominator cannot dominate $G-e$ in two moves.  
Next, let  $e=xa$ (the case $e=xb$ is symmetric).  The $\gamma$-sets of $G-e$ are $\{x,v_1\}, \{x,a\},\{v_1,h_1\}$ and $\{v_1,h_2\}$.  Staller's first move is to select $v_1$.  If now $d_1=x$, then $s_2=a$; if $d_1=a$, then $s_2=x$.  Dominator cannot dominate $G-e$ in one more move.  
Now assume that $e=h_1w$ for some $w\in V(H_m)-\{h_1,h_2\}$. The $\gamma$-sets of $G-e$ are $\{x,v_1\}, \{x,a\},\{x,b\}$ and $\{v_1,h_2\}$.  Staller plays $x$ on her first move.  If $d_1=v_1$, then $s_2=h_2$; otherwise $s_2=v_1$.  Dominator cannot finish dominating $G-e$ in one more move.  A similar strategy will work for Staller if $e=h_2w$ for $w$ as specified above.  (Note that this case does not exist if $m=2$.)
Finally, let $e=az$ for some $z \notin \{v_1,x\}$. The only $\gamma$-sets of $G-e$ are $\{x,v_1\}, \{x,b\},\{v_1,h_1\}$ and $\{v_1,h_2\}$.  Staller will play $v_1$ as her first move.  If $d_1=x$, then $s_2=b$; otherwise $s_2=x$.  A similar strategy works for Staller if $e=bz$ for some $z \notin \{v_1,x\}$.
Therefore, $G$ is $2$-$\gmb'$-critical if $m\neq 3$.

Now let $m=3$ and let $\mathfrak{U}(H_m)=\{h_1,h_2,h_3\}$.  To show that $G$ is critical in this case, the arguments when considering $G-e$ are the same as those above (where $m\neq 3$) with the exception of $e \in \{h_1h_2,h_1h_3,h_2h_3\}$.  Hence, suppose $e=h_1h_2$.  It is easy to show that 
$\{x,v_1\}, \{x,a\},\{x,b\}$ and $\{v_1,h_3\}$ is a complete list of the minimum dominating sets of $G-e$.  Staller can play $x$ for her first move.  If $d_1=v_1$, then $s_2=h_3$; otherwise $s_2=v_1$.  Dominator will then need at least two additional moves to end the game.  The cases $e=h_1h_3$ and $e=h_2h_3$ are symmetric to $e=h_1h_2$.  
Therefore, in all cases we have shown that 
$\gmb'(G-e)>2$ and thus $G$ is $2$-$\gmb'$-critical. 
\end{proof}

\begin{lemma}\label{lem:multipleuniversal} If $G$ is a connected graph with $\gmb'(G)=2$ and has a cut-vertex $x$, then in $G-x = C_1 \cup C_2$ each $C_i$ has a universal vertex in $C_i$ for $i\in[2]$ and at least one component $C_i$ contains at least two universal vertices in $C_i$.
\end{lemma}

\begin{proof}
    If Staller plays $x$ as her first move, then each  $C_i$ contains a universal vertex of $C_i$. Suppose that $C_i$ contains exactly one universal vertex of $C_i$ for $i \in [2]$, say $v_i$ in $C_i$. Without loss of generality we may assume that Dominator first plays in $C_1$. Then Staller's optimal second move is to play $v_2$ in $C_2$, meaning that Dominator cannot dominate $C_2$ in his second move, which is a contradiction.
\end{proof}

The following lemma will be used several times in constructing $2$-$\gmb'$-critical graphs.

\begin{lemma}\label{lem:cut-triangle} Let $G$ be a connected $2$-$\gmb'$-critical graph with a cut-vertex $x$. If $G-x$ contains a component $C$ such that $x$ dominates $C$ and $\mathfrak{u}(C) \ge 2$, then $C = H_m$ for some $m \ge 2$. 
\end{lemma}

\begin{proof}
    Let $G-x = C_1 \cup C_2$ where $x$ dominates $C_2$ and $\mathfrak{u}(C_2) \ge 2$. Let $\{v_2, v_3\} \subseteq \mathfrak{U}(C_2)$. By Lemma~\ref{lem:multipleuniversal}, $\gamma(C_1) = 1$.  Denote a universal vertex of $C_1$ by $v_1$. If $n(C_2) \in \{2, 3\}$, then $C_2 = H_m$ with $m\in \{2, 3\}$. So we shall assume that $n(C_2) \ge 4$. Let $A = V(C_2) - \{v_2, v_3\}$. 

    Since $\gmb'(G)=2$, Dominator has a strategy to win the game in two moves. Hence if Staller first selects $v_1$, then either there is another universal vertex $u \in \mathfrak{U}(C_1) - \{v_1\}$, or Dominator first selects $x$ and there are two more vertices $u_1,u_2 \in \mathfrak{U}(C_1) - \{v_1\}$. In any case, $\mathfrak{u}(C_1) \ge 2$.

    Next, we show that $C_2$ contains exactly two universal vertices. Suppose to the contrary that $C_2$ contains universal vertices $v_2, v_3$, and $w$, and let $y$ be in $V(C_2) - \{v_2, v_3, w\}$. We show that $\gmb'(G - yw) = 2$.  Indeed, if Staller's first move is $v_1$, then Dominator can select $u\in \mathfrak{U}(C_1) - \{v_1\}$ in his first move and finishes the game by playing one vertex from $\{x,v_2,v_3\}$ in his second move. If Staller in her first move does not select $v_1$, then Dominator first plays $v_1$ and he finishes the game by playing one vertex from $\{x,v_2,v_3\}$. Hence, $\gmb'(G-yw)=2$, a contradiction to $G$ being $2$-$\gmb'$-critical. It follows that $C_2$ has exactly two universal vertices.
    
    Finally, we show that $A$ is an independent set. Suppose there is an edge $e$ between two vertices from $A$. Then $\gmb'(G-e)=2$ since Dominator can follow the same strategy as in the case when he was playing the game on $G$ and finishes the game in two moves,  again contradicting the fact that $G$ is $2$-$\gmb'$-critical. Hence, $A$ is an independent set and $C_2=H_m$ where $m=|V(C_2)|$.

\end{proof}

\begin{proposition} 
\label{prp:allbridge}
If $G$ is a connected $2$-$\gmb'$-critical graph with a cut-vertex $x$ where each edge incident to $x$ is a bridge in $G$, then  $G=F_{m_1,1,m_2}$ for some $2 \leq m_1 \leq m_2$.
\end{proposition}

\begin{proof}
By Observation~\ref{o:1}, we may assume that Staller plays vertex $x$ as her first move since every vertex is an optimal first move. Also we know that $\deg_G(x) = 2$ as $\gamma(G-x) = 2$ meaning that $G-x$ contains only two components. Let $G-x = C_1 \cup C_2$. By Lemma~\ref{lem:multipleuniversal} each $C_i$ has a universal vertex in $C_i$. Moreover, since Dominator has exactly two moves which will dominate all of $G$, one move will be in $C_1$ and the other will be in $C_2$ where either the first or second move is a vertex adjacent to $x$. Thus, we may assume, relabeling if necessary, that $v_1$ is a universal vertex in $C_1$ adjacent to $x$. Since $v_1$ is a cut-vertex of $G$, we infer by a similar argument that the two components of $G-v_1$, say $C_1' = C_1 - v_1$ and $\widetilde{C_2}=G[C_2 \cup \{x\}]$, each have a universal vertex. In particular, this implies that there exists a universal vertex, say $v_2$, in $C_2$ which is also universal in $\widetilde{C_2}$ (meaning we can assume $v_2$ is adjacent to $x$). Let $C_2' = C_2 - v_2$. 

Now we consider the orders of $C_i'$ for $i\in [2]$. Note that $n(C_i') \ge 2$ as $G$ cannot contain a leaf by Proposition~\ref{Prop:necessary-conditions}. Moreover, if $n(C_i') =2$, then $C_i$ induces a triangle for the same reasoning. Therefore, if $n(C_i') =2$ for $i \in [2]$, $G \in \mathcal{F}$. Now, we may assume, relabeling if necessary, that $n(C_1') \ge 3$. We first show that $C_i'$ contains at least two universal vertices in $C_i'$ for $i \in [2]$. Since $v_1$ is a cut-vertex in $G$, it follows by Lemma~\ref{lem:multipleuniversal} that one of $C_1'$ or $\widetilde{C_2}$ contains at least two universal vertices. However, $xv_2$ is a bridge and thus the only vertex of $C_2 \subseteq \widetilde{C_2}$ adjacent to $x$ is $v_2$. Hence $\widetilde{C_2}$ contains only one universal vertex $v_2$ and thus $C_1'$ must contain two universal vertices. Interchanging the roles of $C_1$ and $C_2$, we also derive that $C_2'$ contains at least two universal vertices. Note that this implies that if $n(C_i')=3$, then $C_i'$ is a triangle and $C_i$ is a $K_4$. Therefore, we shall assume that $n(C_1') \ge 4$. Now we first apply Lemma~\ref{lem:cut-triangle} to cut-vertex $v_1$ and component $C_1'$ of $G-v_1$ and deduce that $C_1=H_{m_1}$ for some $m_1\geq 3$. In the same way we obtain that $C_2'=H_{m_2}$ for some $m_2 \geq 3$. Since $x$ is a cut-vertex of $G$  of degree $2$ and it is the only neighbor of $v_i$ outside $C_i$ for $i \in [2]$, it follows that $G=F_{m_1,1,m_2}$.
\end{proof}

For the next result we introduce the following definition. If $A$ and $B$ are disjoint subsets of the vertex set of a graph, then by $E(A,B)$ we denote the set of edges with one vertex in $A$ and the other in $B$. 

\begin{lemma}\label{lem:partition} If $G$ is a graph and we can partition $V(G) = V_1 \cup V_2$ such that $\mathfrak{u}(G[V_i])\ge 2$ for $i\in [2]$, then $\gmb'(G)\le 2$. Furthermore, if $|E(V_1, V_2)| \ge 1$, then $G$ is not $2$-$\gmb'$-critical.
\end{lemma}
\begin{proof}
Let $V(G) = V_1 \cup V_2$ be a partition of the vertex set of $G$, where $u_1,u_2\in\mathfrak{U}(G[V_1])$ and $v_1,v_2\in\mathfrak{U}(G[V_2])$. Dominator can use a Pairing Strategy by selecting a vertex from the pair $u_1,u_2$ as the response to Staller selecting the other vertex from that pair, and by selecting a vertex in the pair $v_1,v_2$ as the response to Staller selecting the other vertex from that pair.  (If Staller selects a vertex not in $\{u_1,u_2,v_1,v_2\}$, Dominator selects any vertex from a pair whose vertices have not yet been selected.) Clearly, by this strategy Dominator ensures that at most two vertices will be played during the Maker-Breaker domination game on $G$; that is, $\gmb'(G)\le 2$. 

Now, assume that $E(V_1, V_2)\ne\emptyset$, and let $e$ be an edge with one vertex in $V_1$ and the other in $V_2$. The Pairing Strategy described above ensures that Dominator wins the game on $G-e$ in two moves. Thus $\gmb'(G-e)\le 2$, which implies that $G$ is not $2$-$\gmb'$-critical.
\end{proof}

\begin{proposition} 
\label{prp:v1xbridge}
If $G$ is a connected graph with a cut-vertex $x$ such that $G-x = C_1 \cup C_2$, where $v_1$ is a universal vertex in $C_1$,  $v_1x$ is a bridge in $G$, and $\deg_G(x) \ge 3$,  then $G$ is not $2$-$\gmb'$-critical. 
\end{proposition}

\begin{proof}
    Suppose to the contrary that $G$ is $2$-$\gmb'$-critical. Since $v_1x$ is a bridge, $v_1$ is also a cut-vertex. Let $C_1' = C_1 -v_1$ and $\widetilde{C_2}$ be the component of $G- v_1$ containing $C_2$ and $x$. By Lemma~\ref{lem:multipleuniversal} we know that $C_1'$ and $\widetilde{C_2}$ each have at least one universal vertex and $C_1'$ 
    or $\widetilde{C_2}$
    has at least two universal vertices. Suppose first that $\widetilde{C_2}$ contains at least two universal vertices. Then Lemma~\ref{lem:partition} used for $V_1=V(C_1)$, $V_2=V(\widetilde{C_2})$ implies that $G$ is not $2$-$\gmb'$-critical since there exist edges between $V_1$ and $V_2$, a contradiction. Therefore, we shall assume that $\widetilde{C_2}$ contains exactly one universal vertex, call it $v_2$. Furthermore, if $C_2$ contains another universal vertex other than $v_2$, then it is not adjacent to $x$.  Moreover, since one of $C_1'$ or $\widetilde{C_2}$ contains at least two universal vertices, we may assume that $C_1'$ contains at least two universal vertices. 

    As $\deg_G(x) \ge 3$, there exists $y\in C_2$ adjacent to $x$ where $y$ is not a universal vertex in $C_2$. Since any vertex is an optimal first move for Staller in $G$, assume that Staller chooses $v_2$ first. If Dominator's first move is a vertex in $\widetilde{C_2} - \{v_2\}$, say $w$, then Staller can choose $v_1$ for her second move and now Dominator has to dominate at least one vertex in $\widetilde{C_2} - (N[w] \cup \{v_2\})$ as well as a vertex in $C_1'$ in order to dominate $G$, which is a contradiction. Therefore, we may assume that Dominator's first move is in $C_1$ and by Continuation Principle we may assume he selects $v_1$. It follows that there are at least two other universal vertices in $C_2$ other than $v_2$ in order for $\gmb'(G) = 2$, say $v_3$ and $v_4$.  
    
    We now consider $\gmb'(G-xy)$. If Staller's first move is not $v_2$, then Dominator chooses $v_2$ on his first move and regardless of Staller's second move, there is at least one more universal vertex in $C_1$ that Dominator can choose on his second move so that $G-xy$ is dominated. 
    
    Thus, assume that Staller's first move is $v_2$. But now Dominator can choose $v_1$ as his first move and regardless of Staller's second move, Dominator can still dominate $C_2$ with either $v_3$ or $v_4$. Hence, $\gmb'(G-xy) =2$, contradicting the assumption that $G$ is critical. 
 
\end{proof}

\begin{proposition} 
\label{prp:v1xnotbridge}
Let $G$ be a connected graph with a cut-vertex $x$. If $G-x = C_1 \cup C_2$, $v_1 \in \mathfrak{U}(C_1)$, $\mathfrak{u}(C_1) \le \mathfrak{u}(C_2)$, and $v_1x\in E(G)$ is not a bridge, then $G$ is not $2$-$\gmb'$-critical.
\end{proposition}

\begin{proof}

   Suppose to the contrary that $G$ is $2$-$\gmb'$-critical.  By Lemma~\ref{lem:multipleuniversal}, $\mathfrak{u}(C_2)\geq 2$.
   We let $\{v_2, v_3\} \subseteq \mathfrak{U}(C_2)$. Note that $x$ cannot dominate $C_1$ by Lemma~\ref{lem:partition} using $V_1 = V(C_1) \cup \{x\}$. The same lemma also implies that  $v_1$ is the only vertex in $\mathfrak{U}(C_1)$ adjacent to $x$, yet $xv_1$ is not a bridge, thus there exists $z \in V(C_1) - \mathfrak{U}(C_1)$ which is a neighbor of $x$. 
\vskip2mm
\noindent\textbf{Claim 1:} 
If $v \in V(C_2)$, then $v$ is adjacent to $x$.
\vskip2mm
\noindent\textit{Proof} 
   For the purpose of contradiction suppose that $x$ does not dominate $C_2$. We claim that a vertex in $\mathfrak{U}(C_2)$ is adjacent to $x$. Suppose to the contrary that no such vertex exists. 
   Then Staller's strategy is to first select $v_1$. Then Dominator cannot end the game in only two moves as there is not an available universal vertex of $C_1$ adjacent to $x$ nor a universal vertex of $C_2$ adjacent to $x$. Reindexing if necessary, we may assume $v_2$ is adjacent to $x$. 
   
   Suppose first that $v_2$ is the only vertex in $\mathfrak{U}(C_2)$ adjacent to $x$. When Staller's first move is $v_1$, Dominator must play $v_2$ (otherwise he needs at least three moves to dominate $G$) in his first move. Since $\gmb'(G)=2$ this implies that there are two vertices in $\mathfrak{U}(C_1) - \{v_1\}$ so that Dominator can dominate $C_1$ on his second move. Since $\mathfrak{u}(C_1) \le \mathfrak{u}(C_2)$, this implies that $\mathfrak{u}(C_2) \ge 3$.  We consider $\gmb'(G- xz)$. In $G-xz$, regardless of Staller's first move, Dominator can play either $v_1$ or $v_2$ in his first move which dominates $x$ as well as one of $C_1$ or $C_2$. Since $\mathfrak{u}(C_i) \ge 3$ for $i \in [2]$, Dominator has an optimal second move to finish the game. Thus, $\gmb'(G- xz) = 2$, which is a contradiction. 

Therefore, we may assume that both $v_2$ and $v_3$ are adjacent to $x$. If $\mathfrak{U}(C_1) = \{v_1\}$, then Dominator cannot dominate $G$ in two moves when Staller chooses $v_1$ in her first move. Thus, $\mathfrak{u}(C_1) \ge 2$. Now by using $V_1 = V(C_2) \cup \{x\}$ in Lemma~\ref{lem:partition}, we reach a contradiction. Therefore $x$ dominates $C_2$.
\smallqed

 If $\mathfrak{u}(C_1) >1$, then Lemma~\ref{lem:partition} used for $V_1= V(C_2) \cup \{x\}$ implies that $G$ is not $2$-$\gmb'$-critical, a contradiction. So we may assume that $\mathfrak{U}(C_1) = \{v_1\}$. Note by Proposition~\ref{Prop:necessary-conditions} that $x$ cannot dominate $C_1$ as then $\Delta(G) = n(G)-1$. Let $T \subset V(C_1)$ be the set of vertices in $C_1$ that are not adjacent to $x$. As by Observation~\ref{o:1} $v_1$ is an optimal first move for Staller, in this game Dominator must choose $x$ in his first move for otherwise he cannot dominate $C_1$ in one move played in $C_1$. Since $\gmb'(G)=2$ there exist at least two vertices $a,b$ in $V(C_1)- \{v_1\}$ each of which dominates $T$ so that Dominator can finish the game in two moves. Now we consider $\gmb'(G- xv_1)$. In $G- xv_1$, if Staller does not choose $v_1$ in her first move, then Dominator chooses $v_1$ and finishes the game by choosing one of $x, v_2$, or $v_3$. So assume Staller does choose $v_1$ in her first move. So Dominator's first move is $x$ and then Dominator can play one of the vertices in $\{a,b\}$. However, this implies that $\gmb'(G-xv_1) =2$ which is the final contradiction.

\end{proof}

The last case to consider is when there is no universal vertex in $C_1$ adjacent to $x$ and there are two universal vertices in $C_2$, each adjacent to $x$. 

\begin{proposition} 
\label{prp:xadjacentv2v3}
Let $G$ be a connected $2$-$\gmb'$-critical graph with a cut-vertex $x$ where $G-x = C_1 \cup C_2$. If $\{v_2, v_3\} \subseteq \mathfrak{U}(C_2)$ such that $v_2$ and $v_3$ are adjacent to $x$ and no vertex in $\mathfrak{U}(C_1)$ is adjacent to $x$, then $G \in {\mathcal{F}} \cup {\mathcal{F}}'$.
\end{proposition}
\begin{proof}
First note that if $\mathfrak{u}(C_1) \ge 2$, then $G$ is not $2$-$\gmb'$-critical by Lemma~\ref{lem:partition} using $V_1 = \{x\} \cup V(C_2)$. Therefore, we may assume that $\mathfrak{U}(C_1) = \{v_1\}$  and by assumption, $v_1$ is not adjacent to $x$. If Staller in her strategy first selects $v_1$, then Dominator must play $x$ on his first move for otherwise he cannot dominate $C_1$ in just one move played in $C_1$. This in turn implies that when he plays $x$, $x$ must dominate $C_2$ for otherwise he cannot dominate both $C_1$ and $C_2$ in one more move. Thus, we may assume for the remainder of the proof that $x$ dominates $C_2$. As we have assumed that $\{v_2, v_3\} \subseteq \mathfrak{U}(C_2)$, by Lemma~\ref{lem:cut-triangle} it follows that $C_2 = H_m$ for some $m \ge 2$. 

Let $B$ be the set of neighbors of $x$ in $C_1$. If $V(C_1) - B = \{v_1\}$, then $G \in \mathcal{F}$. (Note that if $B$ is not independent, then by deleting an edge $e$ between two vertices of $B$ we get $\gmb'(G-e)=2$, which is a contradiction.) Therefore, we may assume that $V(C_1) \ne B \cup \{v_1\}$ and let $Y = V(C_1) - (B\cup \{v_1\})$. Assume Staller plays $v_1$ on her first move. As stated above, Dominator must play $x$ on his first move. Since he has a strategy to win in his second move, this means that there exist at least two vertices in $V(C_1)-\{v_1\}$, say $a$ and $b$, each of which dominates $Y$. We have three possibilities:
\begin{enumerate}
    \item  $a,b\in\mathfrak{U}(G[Y])$, or
    \item  $a,b\in B$, or
     \item  $a\in \mathfrak{U}(G[Y])$ and $b\in B$.
\end{enumerate}

{\bf Case 1.} Let $a, b \in \mathfrak{U}(G[Y])$. Suppose first that there exists an edge $cy \in E(G)$ where $c \in B$ and $y \in Y$. We consider $\gmb'(G-cy)$. If Staller doesn't play $v_1$ on her first move in $G - cy$, then Dominator plays $v_1$ in his first move and wins on his second move by chooseing either $v_2$, $v_3$ or $x$. If Staller does choose $v_1$, then Dominator chooses $x$ in his first move and then wins in his second move  by choosing one of $a$ or $b$. Thus, $\gmb'(G- cy) = 2$, contradicting that $G$ is critical. Therefore, we may assume that there is no such edge and $v_1$ is a cut-vertex. By Lemma~\ref{lem:cut-triangle}, $G[Y]= H_r$ for some $r \ge 2$.
Finally, if $B$ is not an independent set in $G$, then $G$ is a supergraph of some $G' \in \mathcal{F}$, implying $G$ is not critical. Thus, $B$ is independent and $G \in \mathcal{F}$. 

{\bf Case 2.} Assume that there exist two vertices $a,b\in B$, each adjacent to all vertices of $Y$. Suppose that there is an edge $e$ in $G[Y]$. When the game is played on $G-e$, $s_1\ne v_1$ gives $d_1=v_1$ and the game clearly ends with two moves by Dominator. On the other hand, if $d_1=x$, then in his second move Dominator can choose either $a$ or $b$ and end the game. Thus $\gmb'(G-e)=2$, and so $G$ is not $2$-$\gmb'$-critical, a contradiction. This implies that $G[Y]$ has no edges. By similar arguments we derive that $G[B]$ has no edges. Analogously, an edge $e$ between $B$ and $Y$ that is incident with neither of $a$ or $b$ implies $\gmb'(G-e)=2$, a contradiction. Thus, letting $q=|B|$, and $s=|Y|$, note that $G$ is isomorphic to $F'_{s,q,m}$, and so $G\in {\cal F}'$.

{\bf Case 3.} Suppose there exists $a\in \mathfrak{U}(G[Y])$ and let $b\in B$ be adjacent to all vertices of $Y$. Since Cases 1 and 2 were already examined, we may assume that $\mathfrak{U}(G[Y])=\{a\}$ and that $b$ is the only vertex in $B$ that is adjacent to all vertices of $Y$. In addition, if $|B|=1$, then $b$ is a universal vertex of $C_1$ adjacent to $x$, a contradiction. Thus $|B|\ge 2$. Now, consider the spanning subgraph $G'=G-xb$ of $G$. Note that the graph $F_{|Y|+1,|B|-1,m} \in \cal F$ is a spanning subgraph of $G'$. Again we infer that $G$ is not $2$-$\gmb'$-critical, a contradiction.
\end{proof}

\bigskip 

We are now ready to complete the proof of the main theorem.   

\medskip

\noindent {\bf Proof of Theorem~\ref{thm:cutvertex}.}
By Propositions~\ref{prp:familyF} and~\ref{prp:familyFprime}, we have that every graph in ${\cal F}\cup {\cal F}'$ is $2$-$\gmb'$-critical. 

For the converse, assume that $G$ is a connected $2$-$\gmb'$-critical and has a cut-vertex $x$.  If all edges incident with $x$ are bridges, then Proposition~\ref{prp:allbridge} gives that $G\in {\cal F}$, as desired. Thus we may assume that not all edges incident with $x$ are bridges. Without loss of generality, let $V(G-x)=C_1\cup C_2$ and $\mathfrak{u}(C_1) \le \mathfrak{u}(C_2)$ (note that $G-x$ has exactly two components, as $\gmb'(G)=2$). 
By Lemma~\ref{lem:multipleuniversal}, $C_1$ has at least one universal vertex, say $v_1$, and $C_2$ has at least two universal vertices. Propositions~\ref{prp:v1xbridge} and~\ref{prp:v1xnotbridge} imply that if $v_1x\in E(G)$, then $G$ is not $2$-$\gmb'$-critical, a contradiction.  Therefore, we may assume that $x$ is not adjacent to a universal vertex of $C_1$. 

Now, we claim that at least two universal vertices of $C_2$, say $v_2$ and $v_3$, are adjacent to $x$. Indeed, if there is only one such vertex, say $v_2\in \mathfrak{U}(C_2)\cap N(x)$, then by selecting $s_1=v_2$, Staller prevents Dominator to  play two vertices, one from $\mathfrak{U}(C_1)$ and the other from $\mathfrak{U}(C_2)$, which dominate $x$ as well as both components $C_1$ and $C_2$. 
Now, since there are (at least) two universal vertices in $C_2$, say $v_2$ and $v_3$, that are adjacent to $x$, Proposition~\ref{prp:xadjacentv2v3} implies that $G\in {\cal F}\cup{\cal F}'$.

\section*{Acknowledgements}
B.B. and T.D. acknowledge the financial support of the Slovenian Research and Innovation Agency (research core funding No.\ P1-0297, projects N1-0285, N1-0431, and the bilateral project between Slovenia and the USA entitled ``Domination concepts in graphs, project No. BI-US/24-26-036''). K.K. would like to thank the financial  support  provided by an AMS-Simons Research Enhancement Grant for Primarily Undergraduate Institution Faculty.

\section*{Declaration of interests}
 
The authors declare that they have no conflict of interest. 

\section*{Data availability}
 
Our manuscript has no associated data.

\end{document}